\documentclass[reqno,11pt]{amsart}
\usepackage{amsmath}
\usepackage{amssymb}
\usepackage{amsthm}
\usepackage{enumerate}
\usepackage[mathscr]{eucal}
\usepackage{adjustbox}
\usepackage{float}
\usepackage{longtable}
\usepackage{lipsum}
\usepackage{amssymb}

\def\mod#1{{\ifmmode\text{\rm\ (mod~$#1$)}
\else\discretionary{}{}{\hbox{ }}\rm(mod~$#1$)\fi}}

\usepackage{float}
\restylefloat{table}

\newtheorem{theorem}{Theorem}[section]
\newtheorem{lemma}[theorem]{Lemma}

\theoremstyle{remark}
\newtheorem{rem}{Remark}

\numberwithin{equation}{section}

\def\mod#1{{\ifmmode\text{\rm\ (mod~$#1$)}
\else\discretionary{}{}{\hbox{ }}\rm(mod~$#1$)\fi}}

\usepackage{hyperref}
\usepackage[normalem]{ulem}
\newcommand{\unclear}[1]{\ifmmode\text{[unclear]}\else\textbf{[unclear: #1]}\fi}

\title{On the Diophantine equation $X^2+p^m=p_1^{n_1}p_2^{n_2}$}
\author{Amir Ghadermarzi}

\begin{document}


\begin{abstract}
We consider the Diophantine equation
\[
X^2+p^m=p_1^{n_1}p_2^{n_2},
\]
where $p,p_1,p_2$ are distinct odd primes. Using factorizations in quadratic fields, sharp bounds for linear forms in two $p$-adic logarithms, elementary methods, and known results on exponential Diophantine equations, we obtain explicit bounds for $m,n_1,n_2$, and consequently the finiteness of the number of solutions for fixed $p,p_1,p_2$.
\end{abstract}

\maketitle

\section{Introduction}
There is an extensive literature on studying the gap between squares and higher powers of integers, reflected in Diophantine equations of the form $x^2+D=y^n$, called Lebesgue--Nagell equations. Many works are devoted to bounding the number of solutions and listing all solutions in special cases by a variety of methods, including methods from algebraic number theory, the hypergeometric method, Baker's theory of linear forms in logarithms, the modularity method, and the primitive divisor theorem of Bilu, Hanrot and Voutier (BHV theorem).

Ap\'ery \cite{apery} proved that for any positive integer $D$ and odd prime $p \nmid D$, the equation $x^2+D=p^n$ has at most two solutions. By applying the BHV theorem \cite{BHV}, Bugeaud and Shorey \cite{bug1} find all cases in which the equation $D_1x^2+D_2=a^n$, where $D_1$ and $D_2$ are coprime positive integers relatively prime to $a$, has more than $2^{\omega(a)-1}$ solutions. In fact, the BHV theorem has been an essential tool in the vast literature on determining solutions of the Lebesgue--Nagell equations in many cases under the condition that $y$ is odd. The BHV theorem has mostly been used to bound $n$ by considering primitive divisors and class number considerations. In some cases, combining the BHV theorem with other methods such as the modularity method and the Thue method has completely solved the equation for small values of $D$ and values of $D$ composed of a finite set of small primes \cite{bennett,siksek}. For a survey of results, see \cite{soydan}.

In comparison, much less work is devoted to equations in which the right-hand side is not a power of an integer. Note that in this case the usual application of the BHV theorem is no longer available. The most general result is due to Evertse \cite{evertse}, who used a reduction of solutions to some unit equations and proved that the equation $x^2+D=p_1^{n_1} p_2^{n_2} \cdots p_r^{n_r}$ has at most $3 \times 7^{4r+6}$ solutions. Yamada \cite{yamada} proved that the equation $X^2+D=2^s p_1^{k_1}p_2^{k_2}$, where $s \in \{0,2\}$ and $p_1,p_2$ are odd primes, has at most 63 solutions.

In this note we investigate the Diophantine equation

\begin{equation} \label{main}
X^2+p^m=p_1^{n_1}p_2^{n_2}.
\end{equation}
Here $p$, $p_1$ and $p_2$ are distinct odd primes. We give explicit bounds on all values $n_i$ and $m$ in terms of the primes involved in the equation. For fixed primes $p$, $p_1$ and $p_2$, this immediately leads to the finiteness of the number of solutions. Our result is more in the spirit of Bugeaud \cite{buguad}, who used linear forms in $p$-adic logarithms to prove explicit bounds for the value $n$ in the solutions of the Diophantine equation $x^2-p^m=\pm y^n$. In order to prove the finiteness of the number of solutions, instead of explicitly bounding $m$, he used an effective result of Kotov \cite{kotov} for the largest prime factor of the polynomial $ax^2+by^n$. The use of $p$-adic linear forms in logarithms in exponential Diophantine equations has its origins in the seminal work of Schinzel \cite{schinzel}.

Before stating our main results, we fix some notation.

Let \[ Q= \min  (p_1,p_2)  , \quad  P= \max  (p_1,p_2 ) \quad S=p_1p_2, \]
  \[\widetilde{P}=\max(p,p_1,p_2)  \quad n= \min(n_1,n_2) \quad N= \max(n_1,n_2). \]
With this notation, we have 
 
 \begin{theorem} \label{theorem asli}
 	Let $(m,n_1,n_2)$ be a solution to \eqref{main}, with at least one of $n_1$ or $n_2$ odd. Then 
$$ 	 N< {9822 Q \widetilde{P}^2 (\widetilde{P} +1)} \log^2(5pS) .$$
 	 \end{theorem}
 	
 	If both $n_1$ and $n_2$ are even, stronger bounds are available, as stated in the following two theorems:

\begin{theorem} \label{first theorem}
	Let $(m,n_1,n_2)$ be a solution to \eqref{main}, with $n_1$ and $n_2$ even and $m$ odd. Then $n< 2 \left( \frac{ \log p}{\log Q} +1\right)$ and $ N <  2(p-1) \left(  \frac{\log p} {\log Q} \right)$.
\end{theorem}

\begin{theorem} \label{second theorem}
	Let $(m,n_1,n_2)$ be a solution to \eqref{main}, with $n_1$, $n_2$ and $m$ even. Then 
\begin{enumerate}
		\item $N< 26 C \log^2(72C)$, where $C=\frac{96p(p+1)}{(\log p)^3 \log(Q)}\log ^2  \left(\widetilde P \right).$
	\item If $ m$  is not a power of $2$ and $ 4 \mid m$, then  $n \leq 4 $ and $N< 2 \max(Q-1,p-1)$. 
\end{enumerate}	

\end{theorem}	

The main feature of these results is the relatively small size of the explicit upper bounds for the exponents $n_1,n_2$ and $m$, in comparison with applying general Baker theory to bound the solutions. To prove these bounds, in Section \ref{m section} we use reduction in the class groups of some imaginary quadratic fields to obtain a linear form in two logarithms, which enables us to apply sharp estimates for linear forms in two $p$-adic logarithms in order to bound $m$ in terms of $N$ and the prime factors in the equation. For the case where $n_1$ and $n_2$ are even, we apply some elementary methods, known results on Diophantine equations, and the derived bound on $m$ to bound $N$. In the final section we study the equation when at least one of $n_1$ and $n_2$ is odd. We apply linear forms in two $p$-adic logarithms, this time to bound the $n_i$'s in terms of $\log m$ and the prime factors involved in the equation. This bound, in conjunction with the bounds for $m$, proves Theorem \ref{theorem asli}. The analogous equation $x^2-p^m=p_1^{n_1}p_2^{n_2}$ requires a separate treatment and is not considered here.

\section{Auxiliary lemmas} \label{auxi}

\begin{lemma} (Ljunggren [\cite {lju}], Størmer [\cite{sto}, p. 168] ) \label{stormer}
The only solution to the equation
\[
x^2+1=2y^n
\]
in integers $x>2$, $y>1$ and $n\geq 3$ is given by
\[
(x,y,n)=(239,13,4).
\]
\end{lemma}

\begin{lemma}[Zsigmondy's theorem \cite{zsi} ] \label{zsigmondy}
Let $a>b \geq 1$ be coprime integers. Then for any integer $n>1$ there is a prime $p$ (called a primitive divisor) that divides $a^n+b^n$ and does not divide $a^k+b^k$ for any $1\leq k<n$ with the exception of
\[
n=3 \qquad \text{and}\qquad a=2,\ b=1.
\]
\end{lemma}

\begin{lemma}[Bennett and Levin   [\cite{Bennett 2}, Proposition 1]  \label{bennett}
If $(x,y,n,q)$ is a solution to
\[
\frac{x^n-1}{x-1}=y^q
\]
with $\lvert x \rvert >1$, $\lvert y \rvert >1$, $n>2$ and $q\geq2$, other than
\[
(3,11,5,2),\quad (7,20,4,2),\quad (18,7,3,3),
\]
\[
(-19,7,3,3),
\]
then $q$ is odd, $|x|>10^4$ and $x$ has a prime divisor $p\equiv 1\pmod q$. 
\end{lemma}

\begin{lemma}[Bugeaud and Laurent; cf. [\cite{buguad}, proposition 1] \label{bugeaud}
Let $p$ be a prime number. Let $\alpha_1$ and $\alpha_2$ be two algebraic numbers which are $p$-adic units. Let $f$ be the residual degree of the extension
\[
\mathbb{Q}_p \subset \mathbb{Q}_p(\alpha_1,\alpha_2)
\]
and put
\[
D=\frac{[\mathbb{Q}(\alpha_1,\alpha_2):\mathbb{Q}]}{f}.
\]
Let $b_1$ and $b_2$ be two positive integers and put
\[
\Lambda=\alpha_1^{b_1}-\alpha_2^{b_2}.
\]
Let $A_1$ and $A_2$ be real numbers greater than $1$ such that
\[
\log A_i\geq \max\left\{h(\alpha_i),\frac{\log p}{D}\right\}
\]
and
\[
b'=\frac{b_1}{D\log A_2}+\frac{b_2}{D\log A_1}.
\]
If $\alpha_1$ and $\alpha_2$ are multiplicatively independent, then we have
\[
v_p(\Lambda)\leq
\frac{24p(p^f-1)}{(p-1)(\log p)^4}D^4
\max\left\{\log b'+\log\log p+0.4,\frac{10\log p}{D},5\right\}^2
\log A_1\log A_2.
\]
\end{lemma}

\begin{lemma} (Peth\"{o} and de Weger \cite{petho} ) \label{petho}
	Let $a \geq 0$ and $h \geq 1$. If $b \geq (e^2/h)^h$ and $x \in \mathbb{R}$ is the largest solution of $x=a+b(\log x)^h$, then 
	$$ x < 2^h  \left( a^{1/h}+b^ {1/h} \log(h^hb) \right)^h$$
\end{lemma}

The following lemma is a modification of Lemma 3 of \cite{buguad}.

\begin{lemma} \label{odd cases}
Let $d>1$ be a squarefree integer, and let $k$ be a positive odd integer, coprime to $d$. Denote by $\ell>1$ the fundamental unit of the field $\mathbb Q(\sqrt d)$. If $X,Y$ and $Z$ are positive integer solutions to the equation
\[
X^2-dY^2=\pm k^Z,
\]
with $(X,Y)=1$, then there exist an integer $t$ and integers $a,b$ and  a positive integer $v$ with $a\equiv b \pmod 2$, and $a$ and $b$ even if $d\not\equiv1\pmod4$, and an algebraic integer $\epsilon= \frac{a+b\sqrt{d}}{2}$ such that
\[
X+Y\sqrt d=\ell^{-t}\left(\frac{a+b\sqrt d}{2}\right)^v.
\]

$$ \lvert  \epsilon \bar{\epsilon} \rvert=k^{z_0}, \quad  \ell^{-1} < \lvert \frac{\epsilon}{\bar{\epsilon}} \rvert<\ell , \quad t< \frac{v	}{2} $$

where $z_0$ divides the class number of $\mathbb Q(\sqrt d)$, and $Z=z_0 v$.

\begin{proof}
	Working in the field $K=\mathbb{Q} (\sqrt{d} )$, since  $k$ and $d$ are coprime and $k$ is odd, the ideals $[X-Y \sqrt{d}]$ and $[X+Y \sqrt{d}]$ are relatively prime in $K$. Therefore, we have the equality of ideals
	
	 $$ [X-Y \sqrt {d}] [X+Y \sqrt{d}] =[k]^Z=  \left ( \mathfrak{a} \mathfrak{\bar{a}}   \right) ^Z$$
where $\mathfrak{a}$ is an integral ideal in $K$. Let $z_0$ be the smallest positive integer such that $\mathfrak{a}^{z_0}$ is a principal ideal in $K$, and let $\alpha$ be a generator of it. Then $z_0$ divides the class number of $K$. Since for any integer $s$, $\alpha l^s$ is also a generator of $\mathfrak{a}^{z_0}$, we can choose a generator $\epsilon$ such that $ \ell^{-1} <  \frac{\epsilon}{\bar{\epsilon}} < \ell$. Putting $v= Z/Z_0$, we obtain
 $$ X+Y \sqrt{d}= \eta \epsilon^ v  \quad  X-Y \sqrt{d}= \bar{\eta} \bar{ \epsilon}^ v$$
 where $\eta$ is a unit in $K$. Therefore, we get 
  $$ X+Y \sqrt{d}=  \ell^{-t}  \left(\frac{a+b\sqrt d}{2}\right)^v$$. 
  From \(X+Y\sqrt d>|X-Y\sqrt d|\) and \(\left|\epsilon/\bar\epsilon\right|\le\ell\), we obtain \(t<v/2\). 
\end{proof}

\end{lemma}

\begin{lemma} (Le \cite{le2} and Louboutin \cite{louboutine)}) \label{class number bound}
Let $D>1$ be an odd squarefree integer. Denote by $h_d$ and $R_d$ the class number and regulator of $\mathbb Q(\sqrt D)$. Then
\[
h_d<D^{1/2}
\qquad\text{and}\qquad
R_d<\frac{\sqrt D}{2}\log (5D).
\]
\end{lemma}

\begin{lemma} (\cite{ramare})\label{class number negative bound}
	Let $D>1$ be a positive squarefree integer. Denote by $h_{-d}$ the class number of $\mathbb Q(\sqrt{-D})$. Then
	 $h_{-d}< \frac{\sqrt{d}}{\pi} \log (20D)$
\end{lemma}
The lemma follows from the class number formula for imaginary quadratic fields and Corollary 1 of \cite{ramare}.
\section{Bounding $m$} \label{m section}
 We start by factoring the equation \eqref{main} in the imaginary quadratic field $K=\mathbb{Q}(\sqrt{-p^m})$.  We obtain:
\[
(x+\sqrt{-p^m})(x-\sqrt{-p^m})=p_1^{n_1}p_2^{n_2}.
\]

The primes $p_1$ and $p_2$ split in $K$, and we have
\[
p_1=\mathfrak a_1\widetilde{\mathfrak a}_1,
\qquad
p_2=\mathfrak a_2\widetilde{\mathfrak a}_2,
\]
where $\mathfrak a_1$ and $\widetilde{\mathfrak a}_1$ are prime ideals over $p_1$, and $\mathfrak a_2$ and $\widetilde{\mathfrak a}_2$ are prime ideals over $p_2$ in $K$.

The vectors $[u,v]\in\mathbb Z^2$ for which
\[
\mathfrak a_1^u\mathfrak a_2^v
\]
is a principal ideal (either integral or fractional) form a lattice in $\mathbb Z^2$. We can choose a basis of the form $[g_1,g_2]$ and $[0,g_3]$ for this lattice, where $0<g_1,g_3<h(K)$, $0\leq g_2<g_3<h_K$, and $g_1g_3=[Z^2:L]\leq h_K$.

Let $\alpha$ be a generator of $\mathfrak a_1^{g_1}\mathfrak a_2^{g_2}$ and let $\beta$ be a generator of $\mathfrak a_2^{g_3}$. Then

\[
x-\sqrt{-p^m}=\zeta\alpha^r\beta^s,
\]

for some root of unity $\zeta\in\mathcal O_K^*$.  Note that by symmetry, we may choose the generators so that \(r\) and \(s\) are both nonzero.
Let
\[
\eta_1=\frac{\overline\alpha}{\alpha},
\qquad
\eta_2=\frac{\overline\beta}{\beta}.
\]

 $$ \frac{2{-p^m}}{x-\sqrt{-p^m}}=  \frac{\bar{\zeta}}{\zeta} \eta_1 ^r  \eta_2^{-s} -1 $$

Raising to the order $\tau$ of $\frac{\bar{\zeta}}{\zeta}$ in $K$ and multiplying by the $p$-adic unit $\eta_2^s$, we obtain the linear form

\begin{equation} \label{p form of m}
\Lambda_{m.p}=\eta_1^{\tau r}-\eta_2^{\tau s}, \quad v_p \left(\Lambda_{m.p} \right) \geq \frac m2.
\end{equation}

\begin{rem}  

Note that if $K\neq\mathbb Q(i),\mathbb Q(\sqrt{-3})$ then $\mathcal O_K^*=\{\pm1\}$ and we may take $\tau=1$. In the cases $K=\mathbb Q(i)$ or $\mathbb Q(\sqrt{-3})$ we have $h_K=1$.

\end{rem}

Since the primes $p_1$ and $p_2$ are both split in $K$, by considering the $\mathfrak a_1$- and $\mathfrak a_2$-adic valuations of the equation $\eta_1^u \eta_2^v=1$, we conclude that they are multiplicatively independent. Hence, in order to bound $m$, we apply Lemma \ref{bugeaud} to \eqref{p form of m} with the following parameters:
$$  \alpha_1=\eta_1= \frac{\overline\alpha}{\alpha} , \quad  \alpha_2= \eta_2= \frac{\overline\beta}{\beta}, \quad b_1= \tau r, \quad b_2= \tau  \lvert s \rvert .$$

Let $f$ be the residual degree of the extension
$
\mathbb{Q}_p \subset \mathbb{Q}_p(\alpha_1,\alpha_2)
$
and put
$
D=\frac{[\mathbb{Q}(\alpha_1,\alpha_2):\mathbb{Q}]}{f}.
$

\begin{rem}
By taking $\alpha_2= \eta_2^{-1}$ if necessary, we can assume that $s$ is positive. 
\end{rem}
With these parameters, we obtain
\[
H_1=h(\alpha_1)=\frac12(g_1\log p_1+g_2\log p_2),
\]
\[
H_2=h(\alpha_2)=\frac12(g_3\log p_2)\leq \frac12 g_3\log Q.
\]
Moreover,
\[
h(\alpha_1)\leq \frac{(g_1+g_2)\log  P}{2} \leq \frac{(g_1+g_3-1) P}{2}
<\frac{g_1g_3\log P}{2}<\frac{h_K\log P}{2},
\]
\[
h(\alpha_2)\leq \frac{h_K\log p_2}{2}.
\]

We choose real numbers $A_i$ such that
\[
\log A_1\geq \max\left\{h(\alpha_1),\frac{\log p}{D}\right\}.
\]
Since $1\leq D\leq2$, we can take
\[
\log A_1=\frac{h_K\log\widetilde P}{D},
\qquad
\text{and}\qquad
\log A_2=\frac{h_K\log\widetilde P}{D}.
\]

To bound the $b_i$'s, note that

\[
b_1= \tau r,  \text{ where } r=\frac{n_1}{g_1}, \text {   and   }
\qquad
b_2= \tau s, \text{where } s=\frac{n_2-(g_2/g_1)n_1}{g_3}.
\]

\[
b_1< \tau N,
\qquad
b_2<2 \tau N.
\]
Hence,
\[
b'=\frac{b_1}{D\log A_2}+\frac{b_2}{D\log A_1}
\]
\[
=\frac{ \tau N}{h_K\log\widetilde P}+
\frac{2 \tau N}{h_K\log\widetilde P}
=\frac{3 \tau N}{h_K\log\widetilde P}.
\]

\[
\log b'+\log\log p+0.4
<\log 3 \tau N-\log h_K+0.4.
\]

On the other hand, $fD=2$ and since $p\geq3$, we obtain $(p^f-1)D^2  \leq (p^2-1)$. Hence, from Lemma \ref{bugeaud}, we get
 
\begin{equation} \label{m bound}
m\leq
\frac{48p(p+1)}{(\log p)^4}
h_K^2\log^2\widetilde P
\max\left\{
\log 3 \tau N+0.4,
\frac{10\log p}{D}
\right\}^2.
\end{equation}
Here, if $\tau>1$, then $\tau\leq6$ and $h_K=1$.
For all integers $N$ and odd primes $p \geq 3$, we have 
$$ \frac{p}{\pi^2} (\log (20p))^2 (\log (3N)+0.4)^2 \geq (\log(18N)+0.4)^2.$$ 
Hence, by Lemma \ref{class number negative bound} (and using $h_K=1$ when $K=\mathbb Q(i)$), we obtain
\begin{equation} \label{m second bound}
 m \leq  
  	\frac{48p^2(p+1)}{(\log p)^4 \pi^2}
  	\log^2 (20p) \log^2\widetilde P
  	\max\left\{
  	\log (3 N)+0.4,
  	\frac{10\log p}{D}
  	\right\}^2.  
  	 \end{equation}
Since $p \geq 3$, 
\begin{equation}  \label{m second ound}
	m \leq  
	\frac{67.55 p^2(p+1)}{(\log p)^2}
 \log^2\widetilde P
	\max\left\{
	\log (3N)+0.4,
	\frac{10\log p}{D}
	\right\}^2.  
\end{equation}	
\section{Both $n_1$ and $n_2$ even}

Throughout this section, we consider solutions to the equation \eqref{main} with $n_1$ and $n_2$ even integers. Let $n_1=2n_1'$ and $n_2=2n_2'$. Then, from equation (1), we have

\begin{equation} \label { zoj factor}
\left(x-p_1^{n_1'}p_2^{n_2'}\right)
\left(x+p_1^{n_1'}p_2^{n_2'}\right)=-p^m.
\end{equation}

Since $p$ is an odd integer and $p,p_1$ and $p_2$ are distinct primes, the factors on the left-hand side of equation \eqref{ zoj factor}  are relatively prime. Therefore,
\[
\begin{cases}
x-p_1^{n_1'}p_2^{n_2'}=-1,\\
x+p_1^{n_1'}p_2^{n_2'}=p^m
\end{cases}
\]
Hence
\begin{equation} \label{ zoj sade shode}
2p_1^{n_1'}p_2^{n_2'}=p^m+1.
\end{equation}

By Zsigmondy's theorem \ref{zsigmondy}, if $m>1$, then either $m=r$, where $r$ is an odd prime, $m=2^k$, or $m=2^k r$, where $r$ is an odd prime. 
\subsection {Proof of Theorem \ref{first theorem} }

To prove Theorem \ref{first theorem}, we assume that $m$ is odd. If $m=1$, the result follows immediately, so we may assume $m>1$, which, as mentioned above, means that $m$ is an odd prime.

From the equation \eqref{ zoj sade shode} we get
\[
(p+1)(p^{m-1}-p^{m-2}+\cdots+1)=2p_1^{n_1'}p_2^{n_2'}.
\]

Let $\beta=(p^{m-1}-p^{m-2}+\cdots+1)$ and
$
d=\gcd(p+1,\beta).
$
Then
\[
p \equiv -1 \pmod d \quad\Rightarrow\quad \beta\equiv m\pmod d.
\]
So $d\mid m$, and since $m$ is a prime, $d=1$ or $d=m$ .

First assume $d=1$. Without loss of generality, we can assume 
\begin{equation} \label{cases 1.1}
\begin{cases}
p+1=2p_1^{n_1'},  \\ 
\beta=\dfrac{p^m+1}{p+1}=p_2^{n_2'}.  
\end{cases}
\end{equation}

$p_2$ is a primitive divisor of $p^m+1$ and $2m\mid p_2-1$.
By the first equation in \eqref{cases 1.1}, we get
\begin{equation} \label{ zoj fard 1}
n_1<\frac{2\log p}{\log p_1}<\frac{2\log p}{\log Q}, 
\end{equation}
By the second equation and Lemma \ref{bennett},
\begin{equation} \label{zoj fard 2}
n_2<2(p-1) \qquad \text{and} \qquad p>10^4.
\end{equation}

Next assume that $d=m$ is prime.

It is clear that $m$ is one of the primes $p_1$ or $p_2$. Say $m=p_1$.  Then $Q=p_1$ and $p_2$ is a primitive divisor of $p^m+1$  and $2p_1 \mid p_2-1$.  Moreover, since 
\[
v_{p_1}(p^m+1)=v_{p_1}(p+1)+v_{p_1}(m),
\]
we get
\begin{equation} \label{cases 1.2}
\begin{cases}
p+1=2p_1^{n_1'-1}  \\
p^{p_1-1}-p^{p_1-2}+\cdots+1=p_1p_2^{n'_2}. 
\end{cases}
\end{equation}

From the first equation of \eqref{cases 1.2}, we get

\begin{equation} \label{zoj fard 3}
4 \leq  n_1<2\left\{\frac{\log p}{\log p_1}+1\right\}
< 2\left\{\frac{\log p}{\log Q}+1\right\}
\end{equation}

By the second equation, since $2p_1\leq p_2$,
\[
(2p_1)^{n'_2}< p_1p_2^{n_2'}<p^{p_1-1}< (2p_1)^{(n'_1-1) (p_1-1)}
\]
Hence,
\[
n'_2 <(p_1-1)(n'_1-1),
\]
Since $2p_1 \leq p+1$,

\begin{equation} \label{zoj fard4}
n_2< (p-1) \left( \frac{ \log p}{\log Q} \right)
\end{equation}

Theorem \ref{first theorem} follows from \eqref{ zoj fard 1}, \eqref{zoj fard 2}, \eqref{zoj fard 3} and \eqref{zoj fard4}.

\subsection {Proof of Theorem \ref{second theorem}}

Throughout this subsection, we assume $m$, $n_1$ and $n_2$ are all even. As mentioned above, this means either $m=2^k$ or $m=2^k r$, where $r$ is an odd prime. The equation \eqref{ zoj sade shode} holds, and since $p^m+1=2p_1^{n'_1}p_2^{n'_2}$, we have:

\begin{equation} \label{N karan zoj}
N <2m \frac{ \log p}{\log Q}.
\end{equation}

We apply the bound in \eqref{m bound} to bound $m$.
Since $m$ is even, we have 
$K=\mathbb Q(i)$, so $h_K=1$ and $ \tau \leq4$.

We obtain
\[
m
<
\frac{48p(p+1)}{(\log p)^4}
\log^2\widetilde P
\max\left\{\log(12N)+0.4,10\log p\right\}^2.
\]
and
\[
N<\frac{96p(p+1)}{(\log p)^3 \log (Q)}\log^2\widetilde P 
\max\left\{\log(12N)+0.4,10\log p\right\}^2.
\]

Let
\[
C=\frac{96p(p+1)}{(\log p)^3 \log(Q)}\log\widetilde P^{\,2}.
\]

If $10\log p > \log(12N)+0.4$, then we get
\begin {equation} \label {first bound for even N}
N< 100 C \log^2 p= 
\frac{9600p(p+1)}{\log (p)\log(Q)}\log\widetilde P^{\,2}
\end{equation}

If
\[
\log (12N)+0.4 > 10 \log p 
\]
then
\[
N<C  \left \{ \log (12N)+0.4 \right\} ^2 < C \log^2 (\exp(0.4).12N)
\]

Therefore, by Lemma \ref{petho}, 
\begin{equation}  \label{ second bound for even N}
 N < 4C \log^2 (72C)
\end{equation}

A direct calculation using the value of $C$ mentioned above and $p \geq 3$  shows that $100 \log(p)^2<26 \log^2(72C)$. Therefore, Theorem \ref{second theorem} follows by \eqref{first bound for even N} and \eqref{ second bound for even N}.

To prove the second part of Theorem \ref{second theorem}, we assume $m=2^k r$ where $k>1$ and $r$ is an odd prime. Let $M=p^{2^k}$. Then by \eqref{ zoj sade shode}  we have
\[
(M+1)(M^{r-1}-M^{r-2}+\cdots+1)=2p_1^{n_1'}p_2^{n_2'}.
\]
Let $\beta= (M^{r-1}-M^{r-2}+\cdots+1)$ and 
$
d=\gcd(M+1,\beta).
$
Then $d=1$ or $d=r$.

 First assume that $d=1$. Then
\begin{equation} \label{m=2^k.r case 1}
\begin{cases}
p^{2^k}+1=2p_1^{n_1'},\\
M^{r-1}-M^{r-2}+\cdots+1=p_2^{n_2'}.
\end{cases}
\end{equation}

By Lemma \ref{stormer},  $n'_1<3$. If $n'_1=2$, since $k\geq2$, we obtain the Diophantine equation $x^4+1=2y^2$ which is known to have no solutions with $x$ or $y$ bigger than $1$. Therefore
\begin{equation} \label{n' value even first}
 n'_1=1 \quad  \text{and} \quad  n_1=2.
\end{equation}
By Lemma \ref{bennett}, we have
\begin{equation} \label{n' value even second}
n'_2<p-1.
\end{equation}
Moreover, $n'_2$ is odd.

Next assume $d\neq1$. Then $d$ is one of the primes $p_1$ or $p_2$, say $p_1$. Then
\begin{equation} \label {2^km second case}
\begin{cases}
M+1=2p_1^{n_1'-1},\\
M^{p_1-1}-M^{p_1-2}+\cdots+1=p_2^{n_2'}.
\end{cases}
\end{equation}
With the same argument as above,

\begin{equation} \label{n' value even third}
n_1'-1=1  \quad\Rightarrow\quad n_1 =4.
\end{equation}
Therefore, the second equation of \eqref {m=2^k.r case 1}  becomes 

$$ \left (2p_1-1 \right)^{p_1-1}-  \left (2p_1-1 \right)^{p_1-2} + \cdots = p_2^{n'_2}$$

But since $p_2$ is a primitive divisor of $p^m+1$, we have $2^{k+1}p_1\mid p_2-1$; therefore,

\begin{equation} \label{n' value even fourth}
M^{n_2'}<M^{p_1-1} \quad\Rightarrow\quad n_2'<p_1-1 
\end{equation}

The second part of Theorem \ref{second theorem} follows from \eqref{n' value even first}, \eqref{n' value even second}, \eqref{n' value even third} and \eqref{n' value even fourth}.

\section{At least one of $n_1$ or $n_2$ is odd}

In this section, we consider \eqref{main} with at least one of $n_1$ or $n_2$ odd and complete the proof of Theorem \ref{theorem asli}. The idea is to bound $n_i$ in terms of $p$, the $p_i$'s and $m$ and then use the bound on $m$ in \eqref{m bound} to obtain explicit bounds on $n_i$. If we specify the parities of the $n_i$'s, we might be able to get stronger results, but we do not pursue that here.

Let

$n_i=2n_i'+r_i ,$ $r_i\in\{0,1\},
$
and
$
d=p_1^{r_1}p_2^{r_2}>1
$
is squarefree. The equation \eqref{main} factors in $\mathbb Q(\sqrt d)$ as
\[
\left(X+p_1^{n_1'}p_2^{n_2'}\sqrt d\right)
\left(X-p_1^{n_1'}p_2^{n_2'}\sqrt d\right)=-p^m.
\]
By Lemma \ref{odd cases}, there exist an algebraic number $\varepsilon=\frac{a+b\sqrt d}{2}$, a positive integer $v$, and an integer $t< \frac{v}{2} $ such that
\begin{equation} \label{cases odd}
\begin{cases}
X+p_1^{n_1'}p_2^{n_2'}\sqrt d=\varepsilon^v\ell^{-t},\\
X-p_1^{n_1'}p_2^{n_2'}\sqrt d= N(\ell)^t\overline\varepsilon^{\,v}\ell^t.
\end{cases}
\end{equation}
Here $\ell$ is the fundamental unit of $\mathbb Q(\sqrt d)$. Putting

\[
\alpha=\frac{\varepsilon}{\bar\varepsilon},
\qquad
\beta=N(\ell)\ell^2.
\]
We obtain
\begin{equation} \label{conjugate}
		\frac{X+Y\sqrt d}{X-Y\sqrt d}
	=\alpha^v\beta^{-t}.
\end{equation}
First assume $t=0$. Then
  $$ X+Y \sqrt{d}= \epsilon^v.$$
Then
  $$ 2Y=bU_v, \quad  U_v= \frac{\epsilon^v - \bar{\epsilon}^v}{\epsilon- \bar{\epsilon}}.$$
Since $\epsilon+\bar\epsilon=a$ and $\epsilon\bar\epsilon=p^{z_0}$ is relatively prime to $a$, the pair $(\epsilon , \bar{\epsilon})$ is a Lucas pair. Consider the related Lucas sequence $U_v$. By the primitive divisor theorem of Bilu, Hanrot and Voutier \cite{BHV}, $U_v$ has a primitive divisor for $v>30$. Every prime divisor of $U_v$ divides $2Y$, while it cannot divide $d$ (related to $v=1$). Hence, either $d=p_1p_2$ and $v\leq30$, or one of the primes $p_1$ or $p_2$ does not divide $d$; say $p_1$. We obtain
   $$ v \leq \max \{ p_1+1, 30\} \leq \frac{30}{7} \widetilde{P} $$
   Moreover, from the first equation in \eqref{cases odd} and the bound $\log\varepsilon<\frac{z\log p+R_d}{2}$ from Lemma \ref{odd cases}, we obtain
   $$ N \log Q < v \left(  z_0 \log p +R_d     \right)$$
   Since $z_0\leq h_d$, the bounds from Lemma \ref{class number bound} give the crude bound 
   $$ N < 19 \widetilde{P}^2 $$
   which is much stronger than the bound stated in Theorem \ref{theorem asli}. Therefore, throughout we assume $t\neq0$.
\begin{lemma} \label{t value}
	Let $\mathbf{N}=\max\{m,N\}$. Then $\left \lvert {t} \right \rvert < \frac{9}{2} \mathbf{N} \frac{\log P}{\log Q}.$
\end{lemma}
\begin{proof}
	By Lemma \ref{odd cases}, $t<\frac{v}{2}$. From \eqref{odd cases} and \eqref{conjugate}, 
	$$ 2t R_d= v \log  \left \lvert \frac{\varepsilon}{\bar\varepsilon} \right \rvert+ m \log p- 2 \log (X+Y\sqrt{d})$$
	Since $X+Y \sqrt{d} < 2Y \sqrt{d} \rightarrow 2 \log (X+Y\sqrt{d}) < N \log(PQ) + 2 \log(2) $
	
	Moreover, from Lemma \ref{odd cases},$ \log  \left \lvert \frac{\varepsilon}{\bar\varepsilon} \right \rvert > -R_d$. On the other hand, $R_d > \frac{1}{4} \log d \geq \frac{1}{4} \log Q $.  Therefore
	 $$t  > - \frac{v}{2}  -2N \frac{\log \left( PQ \right)}{\log Q} +\frac{2( m \log p -2 \log(2))}{\log Q}$$
	 
	Since $m\geq1$ and $p\geq3$,
	 $$ t > - \frac{v}{2}  -2N \frac{\log \left( PQ \right)}{\log Q} -\frac{2\log (4/3)}{ \log Q} $$
	 We may assume $N>P$, since otherwise Theorem \ref{theorem asli} holds. The last term is then absorbed in the preceding estimate. Together with $t< v/2$ and $ v \leq m <\mathbf{N}$ this yields:
	 $$ \left \lvert t \right \rvert < \mathbf{N} \left( \frac{5}{2}+ 2 \frac{\log P}{\log Q} \right) <  \frac{9}{2} \mathbf{N} \frac{\log P}{\log Q}. $$
\end{proof}

From \eqref{cases odd} and \eqref{conjugate}, we obtain $
v_{p_i}\!\left(\alpha^v\beta^{-t}-1\right)=\frac{n_i}{2},
$ and since $\beta$ is a $p_i$-adic unit, this leads to the following linear form in $p_i$-adic logarithms:

\[
\Lambda_{p_i}=\left(\frac{\varepsilon}{\overline\varepsilon}\right)^v-(\beta)^t, \quad 
v_{p_i} \left( \Lambda_{p_i}  \right)= \frac{n_i}{2}
\]

If $t$ is negative, we can replace $\ell$ with $\ell^{-1}$ and let $t$ be positive.

Since $\varepsilon/\overline\varepsilon$ is not a unit, $\varepsilon/\overline\varepsilon$ and $\ell^2$ are multiplicatively independent algebraic numbers that are both $p_i$-adic units. In order to bound $n_i$, we apply Lemma \ref{bugeaud} with
\[
\alpha_1=\frac{\varepsilon}{\overline\varepsilon},\qquad
\alpha_2=\beta,\qquad b_1=v,\qquad b_2=t.
\]
Let $f_i$ be the residual degree of the extension
$
\mathbb{Q}_{p_i }\subset \mathbb{Q}_{p_i}(\alpha_1,\alpha_2)
$
and put
$
D_i=\frac{[\mathbb{Q}(\alpha_1,\alpha_2):\mathbb{Q}]}{f_i}.
$

By Lemma \ref{odd cases}, 
\begin{equation}  \label{epsilon bounds}
h\left(\frac{\epsilon}{\bar\epsilon}\right)
\le
\frac12\left(z_0\log p+
\left|\log\left|\frac{\epsilon}{\bar\epsilon}\right|\right|\right)
<
\frac{z_0\log p+R_d}{2},
\qquad
h\bigl(N(\ell)\ell^2\bigr)=R_d.
\end{equation}

Define
\[
R=\frac{\sqrt{p_1p_2}}{2}\log(5p_1p_2)=\frac{\sqrt S}{2}\log(5S).
\]
Then, by \eqref{epsilon bounds} and Lemma \ref{class number bound}, we can take real numbers $A_1,A_2$ such that
\[
\log A_1=\frac{z_0\log p+R}{2}
\qquad\text{and}\qquad
\log A_2=R.
\]
Since $z_0\log p>0$,
 $$ b'_i < \frac{v+2 \left \lvert t \right \rvert}{D_iR}$$
Therefore, 
 $$\log b'_i + \log \log p_i +0.4 < \log \left( 2 (v+ 2 \left \lvert  t  \right \rvert ) \right) - \log D_i  -\frac{\log P}{2} - \frac{\log Q}{2} + \log \left( \frac{ \log p_i}{\log (5PQ)} \right)+0.4 $$
 By Lemma \ref{t value}, since $D_i\geq1$, $p_i\leq P$ and $v<m\leq\mathbf{N}$, we obtain
 
 $$ \log b'_i + \log \log p_i +0.4 < \log \left( \frac{20  e^{0.4} \mathbf{N} \log P}{\sqrt{PQ} \log Q}  \right )$$
Since $P\geq5$ and $Q\geq3$,
 $$ \log b'_i + \log \log p +0.4 < \log \left( 11.54 \mathbf{N} \right) $$

Therefore,
\[
n_i<
\frac{24p_i(p_i^{f_i}-1)}{(p_i-1)(\log p_i)^4}
D_i^4 R(z\log p+R)
\max\left\{\log  (11.54 \mathbf{N}),\frac{10\log p_i}{D_i}\right\}^2.
\]
Since
\[
z_0 \leq h_d<\sqrt{p_1p_2}
\]
by Lemma \ref{class number bound},
\[
R(z_0\log p+R)
<\frac{\sqrt S}{2}\log(5S)
\left(\sqrt S\log p+\frac{\sqrt S}{2}\log(5S)\right).
\]
Since, for $a,b>1$,
\[
\log a\log b\leq \frac14\log^2(ab),
\]
\begin{equation} \label{ni bound odd}
n_i<
\frac{6S p_i(p_i^{f_i}-1)}{(p_i-1)(\log p_i)^4}
D_i^4\log^2(5pS)
\max\left\{\log  \left(11.54 \mathbf{N} \right) ,\frac{10\log p_i}{D_i}\right\}^2.
\end{equation}

First assume that
\[
\frac{10\log p_i}{D_i}\geq \log \left(11.54 \mathbf{N} \right)
\]
Since $f_iD_i=2$ and $p_i\geq3$,
\[
(p_i^{f_i}-1)D_i^2\leq(p_i^2-1).
\]
Thus, in this case,
\begin{equation} \label{first ni bound}
n_i<
\frac{600S p_i(p_i+1)}{(\log p_i)^2}\log^2(5pS).
\end{equation}

Next assume that
\[
\log  \left(11.54 \mathbf{N} \right)>\frac{10\log p_i}{D_i}.
\]
As $D_if_i=2$,
\[
(p_i^{f_i}-1)D_i^4<4(p_i^2-1).
\]
Let
\[
C_i=\frac{24 S p_i(p_i+1)}{(\log p_i)^4}\log^2(5pS).
\]
Then, in this case,
\begin{equation} \label{Ci va ni}
n_i<C_i\log^2(11.54 \mathbf{N}).
\end{equation}

First assume $N<m$. Then $\mathbf{N}=m$, and from \eqref{m second ound} we have
 $$ N < \frac{67.55 p^2(p+1)}{(\log p)^2}
 ) \log^2\widetilde P
 \max\left\{
 \log 3  N+0.4,
 \frac{10\log p}{D}
 \right\}^2. $$
 
This yields either
 \begin{equation} \label{ni second bound}
 	N< 6755 p^2 (p+1) \log^2 \widetilde{P} . 
 \end{equation}
or
   $$ N< \frac{67.55 p^2(p+1)}{(\log p)^2 }
    \log^2\widetilde P
   \left\{
   \log (3e^{0.4} N)
   \right\}^2.
$$

Since $\frac{p^2(p+1)}{\log p}$ is increasing for $p\geq3$,  $\frac{67.55 p^2(p+1)}{(\log p)^2 }
\log^2\widetilde P <  90.07 \widetilde P^3 $.
Therefore, by Lemma \ref{petho},  
\begin{equation} \label {ni third bound}
N <\frac{270.2 p^2(p+1)}{(\log p)^2 }
\log^2\widetilde P   \log ^2(1613  \widetilde P^3 )
\end{equation}
Finally, assume $\mathbf{N}=N$. It is easy to see that $C_i<102\widetilde{P}^4$. Therefore, by \eqref{Ci va ni} and Lemma \ref{petho}, we obtain
\begin{equation} \label{ni's forth bound}
N< \frac{96 Sp_i(p_i+1)}{\left( \log p_i \right)^4} \log^2(5pS) \log^2(4709 \widetilde{P}^4)
\end{equation} 
Comparing the bounds in \eqref{first ni bound}, \eqref{ni second bound}, \eqref{ni third bound} and \eqref{ni's forth bound}, we obtain
$$ N< \frac{141 Q \widetilde{P}^2 (\widetilde{P} +1)}{\left( \log \widetilde{P} \right)^2} \log^2(5pS) \log^2(4709 \widetilde{P}^4) $$
Since $\widetilde{P}\geq7$,
 $$ N< {9822 Q \widetilde{P}^2 (\widetilde{P} +1)} \log^2(5pS) $$
This completes the proof of Theorem \ref{theorem asli}.

\section*{Acknowledgements}
The author used ChatGPT (OpenAI) for numerical computations and for checking the grammar and wording of the manuscript. The author takes full responsibility for all mathematical arguments and statements.

\end{document}